\documentclass[arxiv,reqno,twoside,a4paper,12pt]{amsart}

\usepackage{tikz}
\usepackage{amsmath, verbatim}
\usepackage{amssymb,amsfonts,mathrsfs}
\usepackage[colorlinks=true,linkcolor=blue,citecolor=blue]{hyperref}
\usepackage[driver=pdftex,heightrounded=true,centering]{geometry}
\usepackage{longtable}
\usepackage{tabularx}
\usepackage{multirow}
\usepackage{caption}

\usepackage[scaled]{helvet} 
\usepackage{courier} 
\usepackage[mathbf]{euler}

\theoremstyle{plain}
\newtheorem{theorem}{Theorem}[section]
\newtheorem{thm}{Theorem}[section]
\newtheorem{prop}[thm]{Proposition}

\newtheorem{defn}[thm]{Definition}

\theoremstyle{definition}

\theoremstyle{remark}
\newtheorem{rem}[thm]{Remark}

         \newtheorem{remark}[theorem]{Remark}

\theoremstyle{plain}

\newcommand{\R}{\mathbb{R}}

\newcommand{\N}{\mathbb{N}}

\newcommand{\supp}{\mathrm{supp}}

\newcommand{\scal}{\mathrm{scal}}
\newcommand{\ric}{\mathrm{Ric}}
\newcommand{\trace}{\mathrm{tr}}

\newcommand{\dv}{\text{ }dV}

\DeclareMathOperator{\cF}{\mathscr{F}}

\DeclareMathOperator{\cH}{\mathscr{H}}
\DeclareMathOperator{\A}{\alpha}
\DeclareMathOperator{\w}{\omega}
\DeclareMathOperator{\V}{\mathcal{V}}
\DeclareMathOperator{\cC}{\mathscr{C}}

\DeclareMathOperator{\ho}{\mathcal{C}^{\A}_{\textup{ie}}}
\DeclareMathOperator{\hho}{\mathcal{C}^{2+\A}_{\textup{ie}}}

\DeclareMathOperator{\hok}{\mathcal{C}^{k, \A}_{\textup{ie}}}

\numberwithin{equation}{section}

\definecolor{qqwuqq}{rgb}{0,0,0}

\begin{document}

\date{\today}

\title[A Survey on the Ricci flow on Singular Spaces]
{A Survey on the Ricci flow on Singular Spaces}

\author{Klaus Kr\"oncke}
\address{University Hamburg, Germany} 
\email{klaus.kroencke@uni-hamburg.de} 

\author{Boris Vertman} 
\address{Universit\"at Oldenburg, Germany} 
\email{boris.vertman@uni-oldenburg.de}

\thanks{Partial support by DFG Priority Programme "Geometry at Infinity"}


\begin{abstract}
In this survey we provide an overview of our recent results concerning 
Ricci de Turck flow on spaces with isolated conical singularities. The
crucial characteristic of the flow is that it preserves the conical singularity. 
Under certain conditions, Ricci flat metrics with isolated conical singularities are stable and  positive scalar curvature is preserved under the flow. We also discuss
the relation to Perelman's entropies in the singular setting, and outline
open questions and future reseach directions.
\end{abstract}

\maketitle

\tableofcontents

\section{Introduction and geometric preliminaries}

Geometric flows, among them most notably the Ricci flow, provide a powerful tool 
to attack classification problems in differential geometry and construct Riemannian 
metrics with prescribed curvature conditions. The interest in this research area only grew
since the Ricci flow was used decisively in the Perelman's proof of Thurston's geometrization 
and the Poincare conjectures. \medskip

The present article summarizes recent results of a continuation of a research program on the 
Ricci flow in the setting of singular spaces, obtained in the papers \cite{Ver-Ricci,KV1,KV2}. The two-dimensional Ricci flow
reduces to a scalar equation and has been studied on surfaces with conical 
singularities by Mazzeo, Rubinstein and Sesum in \cite{MRS} and Yin \cite{Yin:RFO}. 
The Ricci flow in two dimensions is equivalent to the Yamabe flow,
which has been studied in general dimension on spaces with edge singularities
by Bahuaud and the second named author in \cite{BV} and \cite{BV2}. 
\medskip 

In the setting of K\"ahler manifolds, K\"ahler-Ricci flow reduces to a scalar Monge 
Ampere equation and has been studied in case of edge singularities 
in connection to the recent resolution of the Calabi-Yau conjecture on Fano 
manifolds by Donaldson \cite{Donaldson} and Tian \cite{Tian}, see 
also Jeffres, Mazzeo and Rubinstein \cite{JMR}. K\"ahler-Ricci flow in case of isolated conical
singularities is geometrically, though not analytically, more intricate than 
edge singularities and has been addressed by Chen and Wang \cite{Wang1}, 
Wang \cite{Wang2}, as well as Liu and Zhang \cite{LZ}. \medskip

We should point out that in the singular setting, Ricci flow loses its uniqueness
and need not preserve the given singularity structure. In fact, Giesen and 
Topping \cite{Topping, Topping2} constructed a solution to the Ricci flow on 
surfaces with singularities, which becomes instantaneously complete. Alternatively,
Simon \cite{MS} constructed Ricci flow in dimension two and three that smoothens 
out the singularity. \medskip

\emph{Acknowledgements:} The authors thank 
the Geometry at Infinity Priority program of the German Research Foundation 
DFG as well as the Australian Mathematical society 
for its financial support and for providing a platform for joint research. 

\subsection{Isolated conical singularities}

Let us first put up an illustration of a compact Riemannian manifold
with an isolated conical singularity in Figure 1, and then proceed with a precise definition.

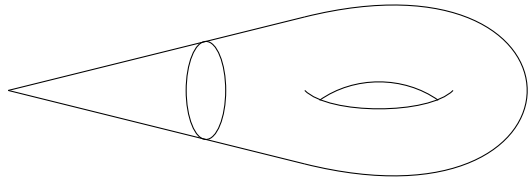
\begin{figure}[h]

\begin{tikzpicture}[scale=1.3]

\draw (-2,4) -- (1,4.75);
\draw (-2,4) -- (1,3.25);
\draw (0,4) ellipse (0.2cm and 0.5cm);
\draw (1,4.75) .. controls (4,5.5) and (4,2.5) .. (1,3.25);
\draw (1,4) .. controls (1.25,3.75) and (2.25,3.75) .. (2.5,4);
\draw (1.15,3.9) .. controls (1.5,4.15) and (2,4.15) .. (2.35,3.9);

\end{tikzpicture}
\label{cone}
\caption{A compact manifold with an isolated conical singularity.}
\end{figure}

\begin{defn}\label{cone-metric}
Consider a compact smooth manifold $\overline{M}$ with boundary 
$\partial M = F$ and open interior denoted by $M$. 
Let $\overline{\cC(F)}$ be a tubular neighborhood of the boundary, 
with open interior $\cC(F) = (0,1)_x \times F$, where $x$ is a 
defining function of the boundary. Consider a smooth Riemannian metric
$g_F$ on the boundary $F$ with $n=\dim F$. An incomplete Riemannian metric $g$ on $M$ 
with an isolated conical singularity is then defined to be smooth away from the
boundary and 
\begin{equation*}
g \restriction \cC(F) = dx^2 + x^2 g_F + h,
\end{equation*}
where the higher order term $h$ has the following asymptotics 
at $x=0$. Let $\overline{g} = dx^2 + x^2 g_F$ denote the exact conical part of the metric $g$ over
$\cC(F)$ and $\nabla_{\overline{g}}$ the corresponding Levi Civita connection. Then we require that
for some $\gamma > 0$ and all integer $k \in \N_0$ the pointwise norm 
\begin{align}\label{higher-order}
| \, x^k \nabla_{\overline{g}}^k h \, |_{\overline{g}} = O(x^\gamma), \ x\to 0. 
\end{align}
\end{defn}

\begin{remark}
We emphasize here that we do not assume that the higher order term $h$
is smooth up to $x=0$ and do not restrict the order $\gamma>0$ to be integer. In that sense
the notion of conical singularities in the present discussion is more 
general than the classical notion of conical singularities where $h$ is usually assumed
to be smooth up to $x=0$ with $\gamma = 1$. This minor generalization is necessary, since the 
Ricci de Turck flow, which will be introduced below, preserves a conical singularity only up to a higher
order term $h$ as above.
\end{remark}

We call $(M,g)$ a compact space with an isolated conical singularity, 
or a \emph{conical manifold} and $g$ a conical metric. The definition naturally extends to conical manifolds 
with finitely many isolated conical singularities.
Since the analytic arguments are local in nature, we may assume without loss
of generality that $M$ has a single conical singularity only.
\medskip

In the present discussion we study compact conical Ricci-flat manifolds $(M,g)$. There are various examples for such spaces. 
Consider a Ricci-flat smooth compact manifold $X$, e.g. a Calabi-Yau manifold or flat torus, 
with a discrete group $G$ acting by isometries, which is not necessarily 
acting strictly discontinuous and admits finitely many fixed points. 
The interior of its quotient $X/G$ defines a compact manifold, an orbifold, with isolated
conical singularities. There exist also examples of compact Ricci-flat manifolds
with non-orbifold isolated conical singularities, constructed by Hein and Sun 
\cite{HeSu}.
\medskip

We now recall elements of b-calculus by Melrose \cite{Mel:TAP, Mel2}.
We choose local coordinates $(x,z)$ on the conical neighborhood 
$\cC(F)$, where $x$ is the defining function of the boundary, $n= \dim F$ and 
$(z)=(z_1,\ldots, z_n)$ are local coordinates on $F$. We
consider the Lie algebra of b-vector fields $\V_b$, which by definition are smooth 
in the interior $M$ and tangent to the boundary $F$. 
In local coordinates $(x,z)$, b-vector fields $\V_b$ are locally generated by 
\[
\left\{x\frac{\partial}{\partial x},  
\partial_z = \left( \frac{\partial}{\partial z_1},\dots, \frac{\partial}{\partial z_n} \right)\right\},
\]
with coefficients being smooth on $\overline{M}$. 
The b-vector fields form a spanning set of section for the 
b-tangent bundle ${}^bTM$, i.e. $\mathcal{V}_b=C^\infty(M,{}^bTM)$. 
The b-cotangent bundle ${}^bT^*M$ is generated locally by the following one-forms
\begin{align}\label{triv}
\left\{\frac{dx}{x}, dz_1,\dots,dz_n\right\}.
\end{align}
These differential forms are singular in the usual sense, but smooth as sections of 
the b-cotangent bundle ${}^bT^*M$. We extend $x:\cC(F) \to (0,1)$ 
smoothly to a non-vanishing function on $M$ and define the incomplete 
b-tangent space ${}^{ib}TM$ by the requirement $C^\infty(M,{}^{ib}TM) = 
x^{-1} C^\infty(M,{}^{b}TM)$. The dual incomplete b-cotangent bundle 
${}^{ib}T^*M$ is related to its complete counterpart by 
\begin{align}
C^\infty(M,{}^{ib}T^*M) = 
x C^\infty(M,{}^{b}T^*M), 
\end{align}
with the spanning sections given locally over $\cC(F)$ by 
\begin{align}\label{triv2}
\left\{dx, x dz_1,\dots, x dz_n\right\}.
\end{align}
With respect to the notation we just introduced, the conical 
metric $g$ in Definition \ref{cone-metric} is a smooth section of 
the vector bundle of the symmetric $2$-tensors of the incomplete
b-cotangent bundle ${}^{ib}T^*M$, i.e. $g \in C^\infty (\textup{Sym}^2({}^{ib}T^*M))$.

\subsection{Ricci de Turck flow and the Lichnerowicz Laplacian}

The Ricci flow is an evolution equation for metrics which reads as
\begin{equation}\label{R}
\partial_t g(t) = -2 \, \textup{Ric} (g(t)).
\end{equation}
Due to diffemorphism invariance of the Ricci tensor, this evolution equation fails to be strongly parabolic. 
One overcomes this problem by adding an additional term to the equation which brakes the diffeomorphism 
invariance. For this reason, one defines the Ricci de Turck flow
\begin{equation}\label{RDT}
\partial_t g(t) = -2 \, \textup{Ric} (g(t)) + \mathcal{L}_{W(g(t),\widetilde{g})} g(t)
\end{equation}
where $W(t)$ is the de Turck vector field defined in terms of the Christoffel symbols for 
the metrics $g(t)$ and a reference metric $\widetilde{g}$ 
\begin{equation}
W(g,\widetilde{g})^k = g^{ij} \left(\Gamma^k_{ij}(g) - \Gamma^k_{ij}(\widetilde{g})\right).
\end{equation}
This flow is equivalent to the Ricci flow via diffeomorphisms.
The linearization of the right hand side is given by 
\begin{equation}
\frac{d}{ds}-2 \, \textup{Ric} (\widetilde{g}+sh) + \mathcal{L}_{W(\widetilde{g}+sh,\widetilde{g})}(\widetilde{g}+sh)|_{s=0}=-\Delta_{L}h,
\end{equation}
where $\Delta_L$ is an elliptic operator which is known as the Lichnerowicz Laplacian. 
Thus, the Ricci de Turck flow is parabolic in the strict sense and, at least in the smooth compact setting, standard existens 
theorems guarantee well-posedness for its initial value problem. Therefore from the analytical perspective, the Ricci de Turck 
flow is much easier to handle than the standard Ricci flow. As it appears in its linearization, the Lichnerowicz Laplacian 
and its spectral properties will be fundamental for considerations in this article. 

\section{Existence of the singular Ricci de Turck flow} 

We shall present here the short time existence result for the singular Ricci de Turck flow obtained by
the second author in \cite{Ver-Ricci} in a simple way, which is sufficient 
for the purpose of the present discussion. 
We consider a compact conical manifold $(M,g_0)$. We study the Ricci de Turck flow with $g_0$ as the 
 initial metric 
While the reference metric $\widetilde{g}$ is usually taken as the initial metric $g_0$, 
in case of $\widetilde{g}$ being Ricci flat, the initial metric $g_0$ can be chosen as a sufficiently small perturbation 
of $\widetilde{g}$.  \medskip

To get shorttime existence for the Ricci de Turck flow
in the conical setting, we need two more conditions. The first one is a condition where we require subquadratic blowup of the Ricci and scalar curvature close to the singular point. The second is a spectral condition on the Lichnerowicz Laplacian which we will explain in the following subsection.

\subsection{Tangential stability}\label{laplace-section}

Let $(M,h)$ be a compact conical Ricci-flat manifold. 
We write $S:= \textup{Sym}^2({}^{ib}T^*M)$. The Lichnerowicz Laplacian 
$\Delta_L: C^{\infty}(M,S)\to C^{\infty}(M,S)$ is a differential operator of second order, 
that can be written in local coordinates near the conical singularity as follows. 
We choose local coordinates $(x,z)$ over the singular neighborhood 
$\cC(F) = (0,1)_x \times F$. In the previous paper \cite{Ver-Ricci} 
we have introduced a decomposition of compactly supported smooth sections $C^\infty_0(\cC(F), 
S \restriction \cC(F))$ 
\begin{equation}
\begin{split}
C^\infty_0(\cC(F), S \restriction \cC(F)) &\to C^\infty_0((0,1), C^\infty(F) 
\times \Omega^1(F) \times \textup{Sym}^2(T^*F)),\\
\w &\mapsto \left( \w(\partial_x, \partial_x), \w (\partial_x, \cdot ), 
\w (\cdot, \cdot)\right),
\end{split}
\end{equation}
where $\Omega^1(F)$ denotes differential $1$-forms on $F$.
Under such a decomposition, the Lichnerowicz Laplace operator 
$\Delta_L$ associated to the singular Riemannian metric $g$ 
attains the following form over $\cC(F)$
\begin{equation}
\Delta_L = - \frac{\partial^2}{\partial x^2} - \frac{n}{x} \frac{\partial}{\partial x}
+ \frac{\square_L}{x^2} + \mathscr{O},
\end{equation}
where $\square_L$ is a differential operator on $C^\infty(F) 
\times \Omega^1(F) \times \textup{Sym}^2(T^*F)$ and the 
higher order term $\mathscr{O} \in x^{-1} \V_b^2$ is a second order
differential operator with one order higher asymptotic behaviour at $x=0$.

\begin{defn}\label{tangential-stability-def}
Let $(F^n,g_F)$ be a closed Einstein manifold\footnote{If $(M,g_0)$ is a conical manifold satisfying condition (2) in Definition \ref{admissible},
then the cross section $(F,g_F)$ of the cone is automatically Einstein with 
Einstein constant $(n-1)$.}  with Einstein constant $(n-1)$.
Then $(F^n,g_F)$ is called (strictly) tangentially stable if the tangential operator 
of the Lichnerowicz Laplacian on its cone restricted to tracefree tensors is 
non-negative (resp. strictly positive).
\end{defn}

\subsection{The existence result}
The conditions for shorttime existence of the Ricci de Turck flow are subsummarised 
under the notion of admissible metrics. In order to improve readability of the article, we put the definitions of the function spaces to an appendix.
\begin{defn}\label{admissible}
Let $(M,g_0)$ be a compact conical manifold.
Then the conical metric $g_0$ is said to be \emph{admissible}, if it satisfies the following assumptions 
for $\gamma>0$ as in \eqref{higher-order}, some $k \in \N$  and $\alpha \in (0,1)$. 
\begin{enumerate}
\item The cross section $(F,g_F)$ is assumed to be tangentially stable.
\item The Ricci curvature has subquadratic growth at the singularitiy, i.e. $\textup{Ric} (g_0) = O(x^{-2+\gamma})$
as $x\to 0$. More precisely, let $\scal(g_0)$ denote the scalar curvature of $g$ and 
$\textup{Ric}^\circ (g_0)$ the trace-free part of the Ricci curvature tensor. Then we assume
\footnote{In view of Definition \ref{cone-metric} \eqref{higher-order} the condition (3) is satisfied
if the leading exact part $\overline{g} = dx^2+x^2g_F$ of the conical metric $g_0$ with 
$g_0\restriction \mathscr{C}(F)=\overline{g} + h$ is Ricci flat, and the higher order term $h$ 
not only satisfies \eqref{higher-order}, but in particular is an element of 
$\mathcal{C}^{k+3,\alpha}_{\textup{ie}}(M,S)_\gamma$. }
\begin{equation}
\begin{split}
&\scal(g_0) \in x^{-2+\gamma} \mathcal{C}^{k+1,\A}_{\textup{ie}}(M, S_1), \\
&\textup{Ric}^\circ (g_0) \in \mathcal{C}^{k+1,\A}_{\textup{ie}}(M, S_0)_{-2+\gamma}.
\end{split}
\end{equation}
\item For any $X_1, \ldots, X_4 \in C^\infty(\overline{M}, {}^{ib}TM)$
we have for the curvature $(0,4)$-tensor 
$$Rm(g_0)(X_1, X_2, X_3, X_4) \in x^{-2} \mathcal{C}^{k+1,\A}_{\textup{ie}}(M).$$ 
\end{enumerate}
\end{defn}

The main result of \cite[Theorem 4.1]{Ver-Ricci}, see also \cite[Theorem 1.2]{KV1}, is then the following theorem. 
\begin{thm}\label{RF}
Let $(M,g_0)$ be a conical manifold with an admissible metric $g_0$. Let the reference
metric $\widetilde{g}$ be either equal to $g_0$ or an admissible conical Ricci flat metric, in which case 
$g_0$ is assumed to be a sufficiently small perturbation of $\widetilde{g}$ in $\cH^{k+2, \A}_{\gamma} (M, S)$.
\medskip

Then there exists some $T>0$, such that the Ricci de Turck flow \eqref{RDT} 
with reference metric $\widetilde{g}$, starting at $g_0$ admits a solution $g(t), t\in [0,T]$, which is 
an admissible perturbation of $g_0$, i.e. $g(t) \in \cH^{k+2, \A}_{\gamma'} (M, S)$ for each $t$, all $k\in \N$ and some
$\gamma' \in (0,\gamma)$ sufficiently small.
\end{thm}
\begin{rem}
A precise definition of the space $\cH^{k+2, \A}_{\gamma'} (M, S)$ is given in the appendix. Let us give a brief description here. Decompose a symmetric $2$-tensor $h$ into its pure trace part and its tracefree part as $h=\trace{h}\cdot g+h_0$. Then $h\in \cH^{k+2, \A}_{\gamma'} (M, S)$ means that as we approach the singularity, $\nabla^l \, h_0=O(x^{-l+\gamma'})$ for $l\in\left\{0,\ldots,k+2\right\}$, $\nabla^l \trace \, {h}=O(x^{-l+\gamma'})$ for $l\in\left\{1,\ldots,k+2\right\}$ but $\trace \, {h}=O(1)$. The last condition distinguishes $\cH^{k+2, \A}_{\gamma'}$ from a standard weighted H\"{o}lder space and ensures that multiples of the metric are also contained in this space.
\end{rem}

This result is obtained as a consequence of a careful microlocal analysis of the heat kernel for the Lichnerowicz 
Laplacian. The heat kernel asymptotics is then used. to establish mapping properties of the heat operator on the
H\"older spaces $\cH^{k+2, \A}_{\gamma}$. It is here, that tangential stability enters in order to obtain these mapping 
properties. Short time existence of the Ricci de Turck flow in these spaces is then a consequence of a fixed point argument,
which requires the metric to be admissible in the sense above to go through.
\medskip

Let us now explain in what sense the flow preserves the conical singularity.
Given an admissible perturbation $g$ of the conical metric $g_0$, the pointwise trace of $g$ 
with respect to $g_0$, denoted as $\textup{tr}_{g_0} g$ is by definition of admissibility an element of the 
H\"older space $\hok (M, S_1)^b_{\gamma}$, restricting at $x=0$ to a constant function $(\textup{tr}_{g_0} g)(0) = u_0>0$.
Setting $\widetilde{x} := \sqrt{u_0} \cdot x$, the admissible perturbation $g= g_0 + h$
attains the form 
$$g = d\widetilde{x}^2 + \widetilde{x}^2 g_F + \widetilde{h},$$
where $|\widetilde{h}|_g = O(x^\gamma)$ as $x\to 0$. Note that the leading part of the 
admissible perturbation $g$ near the conical singularity differs from the leading part of the 
admissible metric $g_0$ only by scaling.

\subsection{Characterizing tangential stability}

A crucial part of our paper \cite{KV1} is devoted to a detailed discussion of the tangential stability. 
Namely, we prove the following general characterization.

\begin{thm}\label{characterizing_tangential_stability}
	Let $(F,g_F)$, $n\geq 3$ be a compact Einstein manifold with constant $(n-1)$. 
	We write $\Delta_E$ for its Einstein operator, and denote the Laplace Beltrami 
	operator by $\Delta$. Then $(F,g_F)$ is tangentially stable if and only if 
	$\mathrm{Spec}(\Delta_E|_{TT})\geq0$ and $\mathrm{Spec}(\Delta)\setminus 
	\left\{0\right\}\cap (n,2(n+1))=\varnothing$. Similarly, $(M,g)$ is strictly tangentially 
	stable if and only if $\mathrm{Spec}(\Delta_E|_{TT})>0$ and $\mathrm{Spec}
	(\Delta)\setminus \left\{0\right\}\cap [n,2(n+1)]=\varnothing$. 
\end{thm}
Establishing this result amounts a careful anaylsis of the Lichnerowicz Laplacian. 
This analysis heavily relies on a decomposition of symmetric two tensors established 
by the first author in \cite{Kro17} to understand the spectrum of $\Delta_L$ on Ricci-flat cones.
\medskip

Any spherical space form is tangentially stable because the Lichnerowicz Laplacian of its cone is the rough Laplacian since the cone is flat. 
However, the spaces $\mathbb{S}^n$ and $\R\mathbb{P}^n$ are not strictly tangentially stable since $2(n+1)\in\mathrm{Spec}
	(\Delta)$ in both cases. 
This property may also hold for other spherical space forms. In the following theorem, we use Theorem \ref{characterizing_tangential_stability}
and eigenvalue computations in \cite{CH15} to characerize tangential stability of symmetric spaces.
\begin{thm}
	Let $(F^n,g_F)$, $n\geq 2$ be a closed Einstein manifold with constant $(n-1)$, 
	which is a symmetric space of compact type. If it is a simple Lie group $G$, 
	it is strictly tangentially stable if $G$ is one of the following spaces:
	\begin{align}
	\mathrm{Spin}(p)\text{ }(p\geq 6,p\neq 7),\qquad \mathrm{E}_6,
	\qquad\mathrm{E}_7,\qquad\mathrm{E}_8,\qquad \mathrm{F}_4.
	\end{align}
	If the cross section is a rank-$1$ symmetric space of compact type $G/K$, 
	$(M,g)$ is strictly tangentially stable if  $G$ is one of the following real Grasmannians
	\begin{equation}
	\begin{aligned}
	&\frac{\mathrm{SO}(2q+2p+1)}{\mathrm{SO}(2q+1)\times \mathrm{SO}(2p)}\text{ }(p\geq 2,q\geq 1),\qquad
	\frac{\mathrm{SO}(8)}{\mathrm{SO}(5)\times\mathrm{SO}(3)},\\
	&\frac{\mathrm{SO}(2p)}{\mathrm{SO}(p)\times \mathrm{SO}(p)}\text{ }(p\geq 4),\qquad
	\frac{\mathrm{SO}(2p+2)}{\mathrm{SO}(p+2)\times \mathrm{SO}(p)}\text{ }(p\geq 4)\\
	&\frac{\mathrm{SO}(2p)}{\mathrm{SO}(2p-q)\times \mathrm{SO}(q)}\text{ }(p-2\geq q\geq 3),
	\end{aligned}
	\end{equation}
	or one of the following spaces:
	\begin{equation}
	\begin{aligned}
	\mathrm{SU}(2p)/\mathrm{SO}(p)\text{ }(n\geq 6),\qquad
	&\mathrm{E}_6/[\mathrm{Sp}(4)/\left\{\pm I\right\}],\qquad \quad
	\mathrm{E}_6/\mathrm{SU}(2)\cdot \mathrm{SU}(6),\\
	\mathrm{E}_7/[\mathrm{SU}(8)/\left\{\pm I\right\}],\qquad&
	\mathrm{E}_7/\mathrm{SO}(12)\cdot\mathrm{SU}(2),\qquad
	\mathrm{E}_8/\mathrm{SO}(16),\\
	\mathrm{E}_8/\mathrm{E}_7\cdot \mathrm{SU}(2),\qquad&
	\mathrm{F}_4/Sp(3)\cdot\mathrm{SU}(2).
	\end{aligned}
	\end{equation}
	\end{thm}

\section{Stability of the singular Ricci de Turck flow}

Our main result in \cite{KV1} establishes long time existence and convergence of the
Ricci de Turck flow for sufficiently small perturbations of conical Ricci-flat metrics, assuming linear and tangential stability and integrability. More precisely we consider a compact conical Ricci-flat manifold $(M,h_0)$  and $g_0$ a sufficiently small perturbation of $h_0$, not necessarily
Ricci-flat. We study the Ricci de Turck flow with $h_0$ as the reference metric, 
and $g_0$ as the initial metric 
\begin{equation}
\partial_t g(t) = -2 \textup{Ric} (g(t)) + \mathcal{L}_{W(t)} g(t), \quad g(0) = g_0,
\end{equation}
where $W(t)$ is 
\begin{equation}
W(t)^k = g(t)^{ij} \left(\Gamma^k_{ij}(g(t)) - \Gamma^k_{ij}(h_0)\right). 
\end{equation}
%
\begin{defn}\label{linear-stability}
We say that $(M,h_0)$ is linearly stable if the the Lichnerowicz 
Laplacian $\Delta_L$ with domain $C^\infty_0(M,S)$ is non-negative.
\end{defn}
\begin{defn}\label{integrability}
	We say that $(M,h_0)$ is \emph{integrable} if for some $\gamma >0$ 
	there exists a smooth finite-dimensional manifold $\cF\subset \cH^{k, \A}_{\gamma} (M, S)$ such that 
	\begin{enumerate}
	      \item  $T_{h_0} \cF = \ker \Delta_{L, h_0} \subset \cH^{k, \A}_{\gamma} (M, S)$,
		\item all Riemannian metrics $h \in \cF $ are Ricci-flat.
		\end{enumerate}
\end{defn}
Our main result is as follows. 
\begin{thm}
Consider a compact conical Ricci-flat manifold $(M,h_0)$. 
Assume that $(M,h_0)$ satisfies the following three additional assumptions
\begin{enumerate}
\item[(i)] $(M,h_0)$ is tangentially stable in the sense of Definition \ref{tangential-stability-def}, 
\item[(ii)] $(M,h_0)$ is linearly stable in the sense of Definition \ref{linear-stability}, 
\item[(iii)] $(M,h_0)$ is integrable in the sense of Definition \ref{integrability}.
\end{enumerate}
If $h_0$ is not strictly tangentially stable, 
we assume in addition that the singularities are orbifold singularities. 
Then for sufficiently small perturbations $g_0$ of $h_0$, there exists a Ricci de Turck flow, with 
a change of reference metric at discrete times, starting at $g_0$ and converging to a conical Ricci-flat 
metric $h^*$ as $t\to\infty$.
\end{thm}
\begin{rem}
If $(M,h_0)$ is a smooth compact manifold, tantential stability is always satisfied. In this case, the statement coincides with the stability results of compact smooth Ricci-flat manifolds obtained in \cite{Ses06}.
\end{rem}
In contrast to the smooth compact case, we can not work with a priori estimates because the curvature is unbounded. Instead, the mapping properties of the heat kernel of the Lichnerowicz Laplacian on the spaces $\cH^{k, \A}_{\gamma} (M, S)$ play a pivotal role in our proof. 
\medskip

We also study examples of compact conical manifolds
where the integrability condition is satisfied. This includes flat spaces with orbifold singularities
as well as K\"ahler manifolds. More precisely we establish the following results. 

\begin{prop}
	Let $(M,h_0)$ be a flat manifold with an orbifold singularity. Then it is linearly stable and integrable.
\end{prop}
The proof of this result is quite simple. Linear stability follows from the absence of curvature and integrability follows from construction the submanifold explicitly as an affine space over $h_0$ modelled over the space of parallel tensors. 
\begin{thm}
	Let $(M,h_0)$ be a Ricci-flat K\"{a}hler manifold where the cross section 
	is either strictly tangentially stable or a space form. 
	Then $h_0$ is linearly stable and integrable.
\end{thm}
This result has been obtained in the smooth compact case in the eighties (see e.g. \cite{Tia87}). The result in our setting follows from carefully adopting techniques from the compact case and using the analysis of weighted Sobolev spaces.

\section{Perelman's entropies on singular spaces}\label{DW-results}
In this section, we review the results obtained in \cite{KV2}
 on Perelman's entropies on compact conical manifolds.
From now on, let $(M^{m},g)$ be a compact conical manifold and $n=m-1$.
We first introduce the three entropies of interest. 

\subsection{The $\lambda$-functional} The $\lambda$-functional is then defined as
\begin{equation*}
\lambda(g) =\inf \left\{\int_M(\scal (g) \omega^2+4|\nabla \omega|^2_g) \dv_g \mid 
\omega \in H^1_1(M), \omega > 0, \int_M \omega^2 \dv_g=1 \right\}.
\end{equation*}
The corresponding Euler-Lagrange equation is
\begin{align}
4\Delta_g \w_g+\scal (g) \w_g=\lambda(g) \w_g,
\end{align}

\subsection{The Ricci shrinker entropy} Consider the functional $\mathcal{W}_-(g,f,\tau)$
\begin{equation*}
\mathcal{W}_-(g,\w,\tau) := \frac{1}{(4\pi\tau)^{m/2}}\int_M [\tau(\scal (g)\cdot \w^2 + 4|\nabla \w|^2_g)
-2\w^2 \ln \w - m \w^2] \dv_g.
\end{equation*}
The Ricci shrinker entropy is then defined by
\begin{equation*}
\mu_-(g, \tau) =\inf \left\{ \mathcal{W}_-(g,\w,\tau) \mid 
\w \in H^1_1(M), \w > 0, \frac{1}{(4\pi\tau)^{m/2}} \int_M \w^2 \dv_g=1 \right\}
\end{equation*}
and the corresponding Euler Lagrange equation is
\begin{align*}
\tau(-4\Delta_g \w_g-\scal (g) \w_g)+2\log(\w_g) \w_g+(m+\nu_-(g, \tau)) \w_g=0.
\end{align*}
It can be shown exactly as in \cite[Corollary 6.34]{CCG}, that if $\lambda(g)>0$, the real number 
$\nu_-(g)=\inf \, \{ \mu_-(g,\tau) \mid \tau>0\}$ exists and is attained by a 
parameter $\tau_g$ and a minimizing function $\w_g$.

\subsection{The Ricci expander entropy}
Consider the functional $\mathcal{W}_+(g,f,\tau)$
\begin{equation*}
\mathcal{W}_+(g,\w,\tau) := \frac{1}{(4\pi\tau)^{m/2}}\int_M [\tau(\scal (g)\cdot \w^2 + \, 4|\nabla \w|^2_g)
+2\w^2 \ln \w + m \w^2] \dv_g.
\end{equation*}
The expander entropy is then defined by
\begin{equation*}
\mu_+(g, \tau) =\inf \left\{ \mathcal{W}_+(g,\w,\tau) \mid 
\w \in H^1_1(M), \w > 0, \frac{1}{(4\pi\tau)^{m/2}} \int_M \w^2 \dv_g=1 \right\}
\end{equation*}
and the corresponding Euler Lagrange equation is
\begin{align}\label{eulerlagrangeexpander}
\tau(-4\Delta_g \w_g-\scal (g) \w_g)-2\log(\w_g) \w_g+(-m+\nu_+(g, \tau)) \w_g=0.
\end{align}
It is now shown exactly as in \cite[p.\ 10]{FIN}, that if $\lambda(g)<0$, the real number $\nu_+(g)=
\sup \, \{\mu_+(g,\tau) \mid \tau>0\}$ exists and is attained by a 
parameter $\tau_g$ and a minimizing function $\w_g$.
\medskip

These functionals are defined almost exactly as in the smooth compact setting. 
The only difference is that one minimizes in the space $H^1_1(M)$ instead of 
$C^{\infty}(M)$. However the weighted Sobolev space $H^1_1(M)$ is exactly 
the right space as it covers the blowup of the scalar curvature at the singular points.
\medskip

The functionals $\lambda$ and $\mu_-$ were already studied by Dai and Wang in the papers \cite{DW1,DW2} in the setting of compact conical manifolds. They showed that they are well-defined and posess minimizers provided that the scalar curvature of the cross section satisfies $\scal(g_F)>m-2$. The minimizers are satisfying for any $\varepsilon > 0$ 
the asymptotics 
\begin{align*}\w_g(x) = o\left(x^{-\frac{m-2}{2} - \varepsilon}\right)\qquad \text{ as } x\to 0.
\end{align*}
By the same methods, we obtained such analogous results for $\mu_+$ under the same assumptions \cite[Theorem 2.16]{KV2}.
\medskip

However, these results can be massively improved by restriction the assumptions on the scalar curvature of the cross section.

\begin{thm}\label{asymptotics_minimizers2}
Let $(M^m,g)$ be a compact conical Riemannian manifold. 
Let $n=m-1$ and $(F^n,g_F)$ be the cross section of the conical part of the metric $g$ and assume that $\scal (g_F)=n(n-1)$.
	Let $\w_g$ be a minimizer in the definition of the $\lambda$-functional, shrinker or the expander entropy. 
	Then there exists an $\overline{\gamma}>0$ such that $\w_g$ admits a partial asymptotic expansion 
	$$\w_g(x,z) = \textup{const} + O(x^{{\overline{\gamma}}}), \ \textup{as} \ x\to 0,$$
	and moreover for $k\in\N$,
	$$|\nabla_g^k\w_g|_g(x,z)=O(x^{{\overline{\gamma}}-k}), \ \textup{as} \ x\to 0.$$
\end{thm}

This result is proved by writing the minimizers in terms of the heat operator, which allows us
to use its mapping properties (Schauder estimates) as already in the proof of the short time
existence of the Ricci de Turck flow. Mapping properties of the heat operator allow us to improve
the asymptotics of the minimizers incrementally and the statement is obtained by an iteration argument.
\medskip

This improvement allows us to study the relation to the Ricci solitons.
Ricci solitions are Riemannian metrics $g$ such that for its 
Ricci curvature $\textup{Ric}(g)$, some vector field $X$, the Lie-derivative 
$\mathcal{L}_X$ and a positive constant $c>0$, the following equations are satisfied
\begin{equation}\label{def_soliton}
\begin{split}
& \textup{Ric}(g) + \mathcal{L}_X g = 0 \quad \textup{(steady Ricci soliton)}, \\
& \textup{Ric}(g) + \mathcal{L}_X g = c\, g \quad \textup{(shrinking Ricci soliton)}.
\end{split}
\end{equation}
Any steady Ricci soliton is up to a diffeomorphism a constant solution to the 
Ricci flow  
\begin{equation}\label{Ricci-steady}
\partial_tg(t) = - 2 \, \textup{Ric} \, (g(t)), \quad g(0)=g.
\end{equation}
Any shrinking Ricci soliton is up to a diffeomorphism a constant solution to the 
normalized Ricci flow
\begin{equation}\label{Ricci-shrinking}
\partial_t g(t) = - 2 \, \textup{Ric} \, (g(t)) + 2c\, g(t), 
\quad g(0)=g.
\end{equation}
Recall also that a Ricci soliton is called gradient if $X=\nabla f$ for some function $f:M\to\R$.
\medskip

Using improved asymptotics we prove the following results, generalizing well known theorems in the compact smooth setting.

\begin{thm}\label{steadysolitonricciflat2}
Let $(M^m,g)$ be a compact conical Riemannian manifold. Let $n=m-1$ and $(F^n,g_F)$ be the cross section of the conical part of the metric $g$. Then the following statements hold. 
\begin{itemize}
\item[(i)] Suppose that $m\geq5$ and $\scal (g_F)=n(n-1)$. Then, if $(M,g)$ is a Ricci soliton, it is gradient. Moreover, if $(M,g)$ is steady or expanding, it is Einstein.
\item[(ii)] In dimension $m=4$, the assertions of part (i) hold if $\ric (g_F)=(n-1)g_F$ .
\end{itemize}
\end{thm}

In addition, we prove monotonicity of the entropies along the singular Ricci de Turck flow.
\begin{thm}\label{monotonicity}Let $(M,g)$ be a compact conical Riemannian manifold of dimension $\dim M \geq 4$ and a tangentially stable cross section.
	\begin{itemize}
		\item[(i)] Then the $\lambda$-functional is nondecreasing along the Ricci de Turck flow preserving conical singularities and constant only along Ricci flat metrics.
		\item[(ii)] Whenever defined, the shrinker and the expanding entropies are nondecreasing along the (normalized) Ricci de Turck flow  preserving conical singularities and constant only along shrinking and expanding solitons, respectively.
	\end{itemize}
	\end{thm}

\section{Positive scalar curvature along singular Ricci de Turck flow}
It is well known that a Ricci flow on smooth compact manifolds preserves the condition of positive scalar curvature. This is an easy consequence of the maximum principle applied to the evolution equation on the scalar curvature. Moreover, the strong maximum principle implies that in this setting a metric of nonnegative scalar curvature that is not Ricci-flat will evolve to a metric of positive scalar curvature immediately. In all these cases, the Ricci flow becomes extinct after finite time. 
\medskip

In current work in progress \cite{KMV}, we establish such results also for the singular Ricci flow on compact conical manifolds. However, we need to strengthen our tangential stability assumption (which guarantees shorttime existence for the singular Ricci flow) somewhat further in order to obtain such results.
In fact we need to assume that $\square_L>n$ on the orthogonal complement of constant functions.
In that case we can obtain existence of singular Ricci de Turck flow in $\cH^{k, \A}_{\gamma} (M\times [0,T], S)$
with $\gamma > 1$. In this case, the de Turck vector field is sufficiently regular (i.e.\ it goes to zero as we approach the conical singularity). This allows us to prove the desired statement.		

\begin{thm}
The singular Ricci de Turck flow preserves positivity of scalar curvature along the flow, provided that the de Turck vector field is sufficiently regular.
Moreover, if the initial metric has nonnegative scalar curvature and is not Ricci-flat, the Ricci de Turck flow will become extinct after finite time.
\end{thm}

The stronger assumption on the tangential operator can be characterized explicitly in a similar manner as in Theorem 
\ref{characterizing_tangential_stability}. For brevity, we only provide a list of symmetric spaces of compact type that 
satisfy that assumption. 
\begin{thm}
		Amongst the symmetric spaces of compact type, only
		\begin{align*}
		E_8,\qquad\mathrm{E}_7/[\mathrm{SU}(8)/\left\{\pm I\right\}],\qquad\mathrm{E}_8/\mathrm{SO}(16),\qquad\mathrm{E}_8/\mathrm{E}_7\cdot \mathrm{SU}(2)
		\end{align*}
		satisfy the conditions $\square_L>n$.\end{thm}

\section{Open questions and further research directions}

One obvious but intricate future research direction is clearly an extension of the 
analysis to non-isolated cones, so-called wedges, and more generally stratified spaces with iterated 
cone-wedge singularities. Already the existence of the various entropies in the edge setting is an open 
question. \medskip

On the other side it is clearly imperative to weaken the conditions of (strict) tangential stability and integrability
for more general applications. This might require a setup of the Ricci de Turck flow in $L^p$ based Sobolev
spaces instead of H\"older spaces as the authors have done till now. This approach would also allow us to 
study the flow of singular metric without a subquadratic blowup of the Ricci curvature. \medskip

Another question is whether we can descend to the Ricci flow. This depends on whether the de Turck vector field
points to or out of the singularity.

\appendix

\section{Sobolev and H\"older spaces}\label{spaces-section} \medskip

Let $\nabla_g$ denote the corresponding Levi Civita covariant derivative. 
Let the boundary defining function $x:\cC(F) \to (0,1)$ be extended smoothly to $\overline{M}$, 
nowhere vanishing on $M$. We consider the space 
$L^2(M)$ of square-integrable scalar functions with respect to the volume form of $g$.
We define for any $s\in \N$ and any $\delta \in \R$ the weighted Sobolev space $H^s_\delta (M)$ as 
the closure of compactly supported smooth functions $C^\infty_0(M)$ under
\begin{align}\label{Sobolev-norm}
\left\|u\right\|_{H^s_\delta} := \sum_{k=0}^s \| x^{k-\delta} 
  \nabla^k_g u \|_{L^2}.
\end{align}
Note that $L^2(M,E) = H^0_{0}(M,E)$ by construction. 

\begin{remark}
An equivalent norm on the weighted Sobolev space $H^s_\delta (M)$ can be defined 
for any choice of local bases $\{X_1, \ldots, X_m\}$ of $\V_b$ as follows. We omit the 
subscript $g$ from the notation of the Levi Civita covariant derivative and write 
\begin{align}
\left\|u\right\|_{H^s_\delta} = \sum_{k=0}^s \sum_{(j_1, \cdots, j_k)} \| \, x^{-\delta} 
\left( \nabla_{X_{j_1}} \circ \cdots \circ \nabla_{X_{j_k}} \right) u \, \|_{L^2}.
\end{align}
\end{remark}

\begin{defn}\label{hoelder-A}
The H\"older space $\ho(M\times [0,T]), \A\in [0,1),$ consists of functions 
$u(p,t)$ that are continuous on $\overline{M} \times [0,T]$ with finite $\A$-th H\"older 
norm
\begin{align}\label{norm-def}
\|u\|_{\A}:=\|u\|_{\infty} + \sup \left(\frac{|u(p,t)-u(p',t')|}{d_M(p,p')^{\A}+
|t-t'|^{\frac{\A}{2}}}\right) <\infty, 
\end{align}
where the distance function $d_M(p,p')$ between any two points $p,p'\in M$ 
is defined with respect to the conical metric $g$, and in terms of the local coordinates 
$(x,z)$ in the singular neighborhood $\cC(F)$ given equivalently by
\begin{align*}
d_M((x,y,z), (x',y',z'))=\left(|x-x'|^2+(x+x')^2|z-z'|^2\right)^{\frac{1}{2}}.
\end{align*}
The supremum is taken over all $(p,p',t) \in M^2 \times [0,T]$\footnote{Finiteness 
of the H\"older norm $\|u\|_{\A}$ in particular implies that $u$ is continuous on the 
closure $\overline{M}$ up to the edge singularity, and the supremum may be taken 
over $(p,p',t) \in \overline{M}^2 \times [0,T]$. Moreover, as explained in 
\cite{Ver-Ricci} we can assume without loss of generality that 
the tuples $(p,p')$ are always taken from within the same coordinate 
patch of a given atlas.}. 
\end{defn} 

We now extend the notion of H\"older spaces to sections of the  
vector bundle $S=\textup{Sym}^2({}^{ib}T^*M)$ of symmetric $2$-tensors. 

\begin{defn}\label{S-0-hoelder}
Denote the fibrewise inner product on $S$ induced by the Riemannian metric $g$, again by $g$.
The H\"older space $\ho (M\times [0,T], S)$ consists of all sections $\w$ of
$S$ which are continuous on $\overline{M} \times [0,T]$, 
such that for any local orthonormal frame 
$\{s_j\}$ of $S$, the scalar functions $g(\w,s_j)$ are $\ho (M\times [0,T])$.
\medskip

The $\A$-th H\"older norm of $\w$ is defined using a partition of unity
$\{\phi_j\}_{j\in J}$ subordinate to a cover of local trivializations of $S$, with a 
local orthonormal frame $\{s_{jk}\}$ over $\supp (\phi_j)$ for each $j\in J$. We put
\begin{align}\label{partition-hoelder-2}
\|\w\|^{(\phi, s)}_{\A}:=\sum_{j\in J} \sum_{k} \| g(\phi_j \w,s_{jk}) \|_{\A}.
\end{align}
\end{defn}

Norms corresponding to different choices of $(\{\phi_j\}, \{s_{jk}\})$
are equivalent and we may drop the upper index $(\phi, s)$ from notation.
The supremum norm $\|\w\|_\infty$ is defined similarly. All the constructions
naturally extend to sections in the sub-bundles $S_0$ and $S_1$, where
the H\"older spaces for $S_1$ reduce to the usual spaces in Definition \ref{hoelder-A}.
\medskip

We now turn to weighted and higher order H\"older spaces. We extend the boundary
defining function $x:\cC(F) \to (0,1)$ smoothly to a non-vanishing function on $M$.
The weighted H\"older spaces of higher order are now defined as follows. 

\begin{defn}\label{funny-spaces}
\begin{enumerate}
\item The weighted H\"older space for $\gamma \in \R$ is
\begin{align*}
&x^\gamma \ho(M\times [0,T], S) := \{ \, x^\gamma \w \mid \w \in \ho(M\times [0,T], S) \, \}, 
\\ &\textup{with H\"older norm} \ \| x^\gamma \w \|_{\A, \gamma} := \|\w\|_{\A}.
\end{align*}
\item The hybrid weighted H\"older space for $\gamma \in \R$ is
\begin{align*}
&\ho_{, \gamma} (M\times [0,T], S) := x^\gamma \ho(M\times [0,T], S)  \, \cap \, 
x^{\gamma + \A} \mathcal{C}^0_{\textup{ie}}(M\times [0,T], S) \\
&\textup{with H\"older norm} \  \| \w \|'_{\A, \gamma} := \|x^{-\gamma} 
\w\|_{\A} + \|x^{-\gamma-\A} \w\|_\infty.
\end{align*}
\item The weighted higher order H\"older spaces, which specify regularity of solutions 
under application of the Levi Civita covariant derivative $\nabla$ of $g$ on symmetric $2$-tensors
and time differentiation are defined for any 
$\gamma \in \R$ and $k \in \N$ by\footnote{Differentiation is a priori 
understood in the distributional sense.}
\begin{equation*}
\begin{split}
&\hok (M\times [0,T], S)_\gamma = \{\w\in \ho_{,\gamma} \mid  \{\nabla_{\V_b}^j \circ \, (x^2 \partial_t)^\ell\} 
\, \w \in \ho_{,\gamma} \ \textup{for any} \ j+2\ell \leq k \}, \\
&\hok (M\times [0,T], S)^b_\gamma = 
\{u\in \ho \mid  \{\nabla_{\V_b}^j \circ \, (x^2 \partial_t)^\ell\} \, u \in 
x^\gamma\ho \ \textup{for any} \ j+2\ell \leq k\},
\end{split}
\end{equation*}
where the upper index b in the second space indicates the fact that despite
the weight $\gamma$, the solutions $u \in \hok (M\times [0,T], S)^b_\gamma$
are only bounded, i.e. $u\in \ho$. The corresponding H\"older norms are defined 
using local bases $\{X_i\}$ of $\V$ and $\mathscr{D}_k:=\{\nabla_{X_{i_1}} \circ \cdots \circ 
\nabla_{X_{i_j}} \circ (x^2 \partial_t)^\ell \mid j+2\ell \leq k\}$ by
\begin{equation*}
\begin{split}
&\|\w\|_{k+\A, \gamma} = \sum_{j\in J} \sum_{X\in \mathscr{D}_k} \| X (\phi_j \w) \|'_{\A, \gamma}
+ \|\w\|'_{\A, \gamma}, \quad \textup{on} \ \hho (M\times [0,T], S)_\gamma, \\
&\|u\|_{k+\A, \gamma} = \sum_{j\in J} \sum_{X\in \mathscr{D}_k} \| X (\phi_j u) \|_{\A, \gamma}
+ \|u\|_{\A}, \quad \textup{on} \ \hho (M\times [0,T], S)^b_\gamma.
\end{split}
\end{equation*} 
\item In case of $\gamma=0$ we just omit the lower weight index and write
e.g. $\hok (M\times [0,T], S)$ and $\hok (M\times [0,T], S)^b$.
\end{enumerate}
\end{defn} 

The H\"older norms for different choices of local bases $\{X_1, \ldots, X_m\}$ of $\V_b$ and different choices
of conical Riemannian metrics $g$ are equivalent due to compactness of $M$ and $F$. \medskip

The vector bundle $S$ decomposes into a direct sum of sub-bundles
\begin{align}
S= S_0 \oplus S_1, 
\end{align}
where the sub-bundle $S_0=\textup{Sym}_0^2({}^{ib}T^*M)$
is the space of trace-free (with respect to the fixed metric $g$) symmetric $2$-tensors,
and $S_1$ is the space of pure trace (with respect to the fixed metric $g$) symmetric 
$2$-tensors. The sub bundle $S_1$ is trivial real vector bundle over $M$ of rank 1.
\medskip

Definition \ref{funny-spaces} extends ad verbatim to sections of $S_0$ and $S_1$.
Since the sub-bundle $S_1$ is a trivial rank one real vector bundle, its sections
correspond to scalar functions. In this case we may omit $S_1$ from the notation and
simply write e.g. $\hok (M\times [0,T])^b_\gamma$. \medskip

\begin{remark}
The higher order weighted H\"older spaces in Definition \ref{funny-spaces}
differ slightly from the corresponding spaces in \cite{Ver-Ricci} by the choice of 
admissible derivatives. While in \cite{Ver-Ricci} we allow differentiation by any 
b-vector field $\V \in \V_b$, here we employ only derivatives of the form $\nabla_{\V},
\V \in \V_b$. \medskip
\end{remark}

Below we will simplify notation by introducing the following spaces.

\begin{defn}\label{H-space}
Let $(M,g)$ be a compact conical manifold and 
assume that the conical cross section $(F,g_F)$ is
strictly tangentially stable. Then we define
\begin{equation*}
\cH^{k, \A}_{\gamma} (M\times [0,T], S) := \hok (M\times [0,T], S_0)_{\gamma}
\ \oplus \ \hok (M\times [0,T], S_1)^b_{\gamma}.
\end{equation*}
If $(F, g_F)$ is tangentially stable but not 
strictly tangentially stable, we set instead
\begin{equation*}
\cH^{k, \A}_{\gamma} (M\times [0,T], S) := \hok (M\times [0,T], S)^b_{\gamma}.
\end{equation*}
\end{defn}

\end{document}